\ifx\shlhetal\undefinedcontrolsequence\let\shlhetal\relax\fi

\documentclass[11pt]{amsart}

\usepackage{amsmath}
\usepackage{amssymb}

\newtheorem{theorem}{Theorem}[section]
\newtheorem{claim}[theorem]{Claim}

\newtheorem{proposition}[theorem]{Proposition}
\newtheorem{corollary}[theorem]{Corollary}

\theoremstyle{definition}
\newtheorem{definition}[theorem]{Definition}

\newtheorem{question}[theorem]{Question}

\theoremstyle{remark}
\newtheorem{remark}[theorem]{Remark}

\newtheorem{conclusion}[theorem]{Conclusion}

\newcount\skewfactor
\def\mathunderaccent#1#2 {\let\theaccent#1\skewfactor#2
\mathpalette\putaccentunder}
\def\putaccentunder#1#2{\oalign{$#1#2$\crcr\hidewidth
\vbox to.2ex{\hbox{$#1\skew\skewfactor\theaccent{}$}\vss}\hidewidth}}

\def\smallbox#1{\leavevmode\thinspace\hbox{\vrule\vtop{\vbox
   {\hrule\kern1pt\hbox{\vphantom{\tt/}\thinspace{\tt#1}\thinspace}}
   \kern1pt\hrule}\vrule}\thinspace}

\newcommand{\bool}{{\bf B}}

\newcommand{\cf}{{\rm cf}}

\newcommand{\Depth}{{\rm Depth}}

\newcommand{\st}{{such that}}
\newcommand{\seq}{{sequence}}
\newcommand{\cont}{{continuous}}
\newcommand{\incr}{{increasing}}

\newcommand{\Wlog}{{Without loss of generality}}

\newcommand{\then}{{\underline{then}}}
\newcommand{\Then}{{\underline{Then}}}

\def\qedref#1{$\qed_{\reforiginal{#1}}$}

\title{$(\kappa,\theta)$-weak normality}
\author{Shimon Garti and Saharon Shelah}
\address{Institute of Mathematics
 The Hebrew University of Jerusalem
 Jerusalem 91904, Israel
 and  Department of Mathematics
 Rutgers University
 New Brunswick, NJ 08854, USA}
\email{shelah@math.huji.ac.il}
\urladdr{http://www.math.rutgers.edu/\char`\~shelah}
\thanks{First typed: November 2009 \\
Research supported by the United States-Israel Binational
Science Foundation. Publication 956 of the second author}
\subjclass[2000]{Primary: 03E04, Secondary: 03G05}
\keywords{Set theory, ultrafilters, weak normality, 
Boolean Algebras, Depth, measurable cardinal}

\begin{document}
\let\labeloriginal\label
\let\reforiginal\ref

\begin{abstract}
We deal with the property of weak normality (for non-principal ultrafilters). We characterize the situation of $|\prod \limits_{i<\kappa} \lambda_i / D|=\lambda$. We have an application for a question of Depth in Boolean Algebras.
\end{abstract}

\maketitle

\newpage

\section{introduction}
The motivation of this article, emerged out of a 
question about the Depth of Boolean Algebras.
We found that a necessary condition to a positive answer 
on a question of Monk (appears in \cite{MR1393943}) depends
on the following condition. We need a \seq\ of cardinals 
$\bar{\lambda} = \langle \lambda_i : i < \kappa\rangle$ with limit $\lambda$ (or just $\lambda = {\rm lim}_D (\bar \lambda)$, see definition \ref{0.4} below, and for simplicity $i < \kappa \Rightarrow \lambda_i \leq \lambda$), 
and an ultrafilter $D$ on $\cf(\lambda) = \kappa$, \st\ 
$|\prod\limits_{i<\kappa} \lambda_i/D|=\lambda$ 
(see \cite{MR2217966} and \cite{temp911}, about the connection 
to Boolean Algebras; We give new results about the 
Depth, in \S2).

These requirements are purely set-theoretical, and they 
depend on the nature of $\kappa$ and $\lambda$, and also 
on the properties of $D$. On one hand, if $D$ is a regular 
ultrafilter then $|\prod\limits_{i<\kappa}\lambda_i/ D|=
\lambda^\kappa$. Notice that $\lambda^\kappa>\lambda$ 
in our case, since $\cf(\lambda)\leq \kappa$. On 
the other hand, having a measurable cardinal $\kappa = \cf(\lambda)$ (or just $\cf(\lambda) \leq \kappa$, $\kappa$ is measurable) 
and a normal ultrafilter $D$, we can choose a \seq\ 
as above, with $|\prod\limits_{i<\kappa}
\lambda_i/D|=\lambda$.

Regular ultrafilters and normal ultrafilters are two poles. 
The question is, what happens to other creatures in the zoo of ultrafilters. We will introduce here the notion of weak normality (the basic notion appears in \cite{MR506523}, and the general notion is taken from \cite{MR1261217}),
and prove two theorems. First, $|\prod\limits_{i<\kappa}
\lambda_i/ D|=\lambda$ implies that $D$ is weakly 
normal (in the sense of definition \ref{0.3} below). 
Second, that under the assumption of weak normality one can 
find $\bar{\lambda} = \langle \lambda_i:i<\kappa\rangle$ with the properties above.

Recall that a normal ultrafilter on $\kappa$ is closed under diagonal intersections of $\kappa$ sets from the ultrafilter.
It follows, that any regressive function $f$ on $\kappa$,
has a suitable set $S_f$ in the ultrafilter, \st\ $f$ is 
constant on $S_f$. In other words, one can find a 
(unique) ordinal $\alpha_*$, \st\ 
$\{i<\kappa:f(i)=\alpha_*\}\in D$ (when $D$ is the 
normal ultrafilter).

This property of regressive functions, leads us to 
another notion of normality. It might happen that 
for no $\alpha_*$ one can get $\{i<\kappa:f(i)=
\alpha_*\}\in D$, but for some $\alpha_*<\kappa$ we 
have $\{i<\kappa:f(i)\leq \alpha_*\}\in D$.

\begin{definition}
\label{0.1} 
Weak normality. \newline 
Let $\kappa$ be an infinite cardinal, $D$ a uniform 
ultrafilter on $\kappa$. \newline 
We say that $D$ is weakly normal, when:
\begin{enumerate}
\item[(*)] For every regressive function $f$ on $\kappa$, 
one can find $\alpha_*<\kappa$, \st\ 
$\{i<\kappa:f(i)\leq \alpha_*\}\in D$
\end{enumerate}
\end{definition}

Every normal ultrafilter is also weakly 
normal. The opposite need not to be true. 
If $D$ satisfies the weak normality condition of 
$\leq \alpha_*$, but not the requirement of 
$=\alpha_*$, then $D$ is not $\kappa$-complete, 
so it is not a normal ultrafilter. 

For our needs, we would like to generalize the notion 
of weak normality. So far, we focused on regressive 
functions from $\kappa$ into $\kappa$. Let us define the 
property of regressiveness, in a more general context.

\begin{definition}
\label{0.2} 
Regressive pairs. \newline 
Let $(\kappa,\theta)$ be a pair of cardinals, $D$ 
an ultrafilter on $\kappa$. \newline 
Let $g : \kappa \rightarrow \theta$ be any function. \newline
We say that $f:\kappa\rightarrow \theta$ is $(\kappa,g)$-regressive, if $i<\kappa\Rightarrow f(i)<g(i)$.
\end{definition}

In the light of definition \ref{0.1}, taking $\theta=\kappa$
and $g \equiv {\rm id}_\kappa$ gives the familiar notion 
of a regressive function on $\kappa$. We would like to form the 
new concept of weak normality, based on the regressive functions of 
\ref{0.2}. But look, if we choose $g(i)=0$ for any
$i<\kappa$, or even $g:\kappa\rightarrow \theta$ bounded,
then we will have an uninteresting definition.
That's the reason for demand (i) in part (a) below:

\begin{definition}
\label{0.3}
$(\kappa,\theta)$-weak normality. \newline 
Let $(\kappa,\theta)$ be a pair of cardinals, 
$g:\kappa \rightarrow \theta$, and $D$ an ultrafilter
on $\kappa$.
\begin{enumerate}
\item[(a)] $D$ is $(\kappa,g)$-weakly normal, if
\begin{enumerate}
\item[(i)] $\epsilon<\theta \Rightarrow \{i<\kappa:g
(i)\geq \epsilon\}\in D$
\item[(ii)] For any $(\kappa,g)$-regressive 
$f$, there is $j_f<\theta$, \st\ $\{i<\kappa:f(i)
< j_f\}\in D$
\end{enumerate}
\item[(b)] $D$ is $(\kappa,\theta)$-
weakly normal if there is a function $g:\kappa\rightarrow \theta$ 
\st\ $D$ is $(\kappa,g)$-weakly normal.
\end{enumerate}

Two remarks about the definition. First, we speak 
about an ultrafilter (that's what we need for our claims), 
but the definition (with some modifications) applies also to a filter. Second, we use $f(i)<j_f$ (instead of $\leq$ in \ref{0.2}), but there
is no essential difference. 

The last definition that we need, adapts the notion
of limit for sequence of cardinals to the notion of 
an ultrafilter. 
\end{definition}

\begin{definition}
\label{0.4}
${\rm lim}_D (\bar \lambda)$. \newline 
Let $\bar \lambda=\langle \lambda_i:i<\kappa\rangle$ 
be a \seq\ of cardinals, $D$ an ultrafilter on $\kappa$. \newline
$\mu:= {\rm lim}_D (\bar \lambda)$ is the (unique) 
cardinal \st\ $\{i<\kappa:\beta<\lambda_i\leq \mu\}\in D$,
for every $\beta<\mu$.
\end{definition}

\par \noindent We conclude this section with some elementary facts.

\begin{claim}
\label{0.5}
Assume $\bar \mu=\langle \mu_j:j<\theta\rangle$ 
is an \incr\ \seq\ of cardinals, with limit $\lambda$. 
Let $D$ be a $(\kappa,\theta)$-weakly normal ultrafilter 
on $\kappa$, and $g:\kappa\rightarrow \theta$ a witness.
Let $\lambda_i=\mu_{g(i)}$, for every $i<\kappa$. \newline
\Then\ $\lambda={\rm lim}_D 
(\langle \lambda_i:i<\kappa\rangle)$.
\end{claim}

\par \noindent \emph{Proof}. \newline 
Easy, by the definition of ${\rm lim}_D$. \newline 

\hfill \qedref{0.5}

\medskip

In the following claim we learn something about the 
relationship between ${\rm lim}_D (\langle 
\lambda_i:i<\kappa\rangle)$ and 
$|\prod\limits_{i<\kappa}\lambda_i/D|$:

\begin{claim}
\label{0.6}
$|\prod\limits_{i<\kappa} \lambda_i/D|\geq 
{\rm lim}_D (\langle \lambda_i:i<\kappa\rangle)$.
\end{claim}

\par \noindent \emph{Proof}. \newline 
Assume to contradiction, that 
$|\prod\limits_{i<\kappa} \lambda_i/D|=
\mu<{\rm lim}_D (\bar \lambda)$.
Choose $\beta<{\rm lim}_D (\bar \lambda)$ \st\ 
$\mu<\beta$. By \ref{0.4} we have: 

$$
A:=\{i<\kappa:\mu<\beta<\lambda_i\leq {\rm lim}_D 
(\bar \lambda)\}\in D.
$$

Define $\chi={\rm Min}\{\lambda_i:i\in A\}$.
Easily, one can define a \seq\ 
$\langle a_\alpha:\alpha<\chi\rangle$ of 
members in $\prod\limits_{i<\kappa} 
\lambda_i /D$, \st\ $a_\alpha<_D a_\beta$ (notice that one needs to define the $a_\alpha$-s only on the set $A$, and $0$ on $\kappa \setminus A$ is alright). 
But $\chi>\beta>\mu$, contradicting the fact 
that $|\prod\limits_{i<\kappa} \lambda_i/D|
=\mu$. \newline 

\hfill \qedref{0.6}

\medskip

We say that $(\prod\limits_{i<\kappa}
\lambda_i,\leq_D)$ is $\theta$-directed, 
if any $A\subseteq \prod\limits_{i<\kappa} 
\lambda_i/D$ satisfies 
$|A|<\theta\Rightarrow A$ has an upper bound 
in $\prod\limits_{i<\kappa} \lambda_i/D$. We say that $\theta$ is $\kappa$-strong when $\alpha < \theta \Rightarrow |\alpha|^\kappa < \theta$. 
The following useful claim draws a line between 
$\theta$-directness and the cardinality of 
$\prod\limits_{i<\kappa}\lambda_i/D$.

\begin{claim}
\label{0.7}
Simple properties of cardinal products. \newline 
Let $D$ be an ultrafilter on $\kappa$.
\begin{enumerate}
\item[(a)] If $(\prod\limits_{i<\kappa} 
\lambda_i,\leq_D)$ is $\theta$-directed,
\then\ $|\prod\limits_{i<\kappa}
\lambda_i/D|\geq \theta$
\item[(b)] If $\kappa_i=\cf(\kappa_i)$ 
for every $i<\kappa$, and $\beta < \theta \Rightarrow \{ i < \kappa : \beta < \kappa_i \} \in D$, \then\ 
$(\prod\limits_{i<\kappa} \kappa_i,\leq_D)$ 
is $\theta$-directed 
\end{enumerate}
\end{claim}

\par \noindent \emph{Proof}.
\begin{enumerate}
\item[(a)] Easy, since if $|\prod\limits_{i<\kappa}
\lambda_i/D|=\theta_*<\theta$, then there exists an unbounded \seq\ of members in $\prod\limits_{i<\kappa}\lambda_i/D$, of 
length $\theta_*$, contradicting the $\theta$-directness.
\item[(b)] Having $A\subseteq \prod\limits_{i<\kappa}
\kappa_i/D$, $|A|=\theta_*<\theta$, just take 
the supremum of $g(i)$ for every $g\in A$, on the set $\{ i<\kappa : \theta_* < \kappa_i \}$ (and $0$ on the rest of the $i$-s). By our assumptions, we get an upper bound for the set $A$ which
belongs to $\prod\limits_{i<\kappa} \kappa_i/D$.
\end{enumerate}

\hfill \qedref{0.7}

\medskip

The last proposition that we need, is about the 
connection between $\theta=\cf(\lambda)$ and $\lambda$.
We defined the property of $(\kappa,\theta)$-weak 
normality, when $\theta=\cf(\lambda)$. We 
concentrated in $(\kappa,g)$-regressive functions, when
$g:\kappa\rightarrow \theta$.
But sometimes we want to pass from $\theta$ to $\lambda$
in our treatment.

\begin{claim}
\label{0.8}
Let $\theta=\cf(\lambda)\leq\kappa<\lambda,D$
an ultrafilter on $\kappa,
\langle \mu_j:j<\theta\rangle$ \incr\ continuous with 
limit $\lambda,g:\kappa\rightarrow \theta$ and 
$\lambda_i=\mu_{g(i)}$ for every $i<\kappa$
\st\ ${\rm lim}_D (\langle \lambda_i:i<\kappa\rangle)
=\lambda$.
Assume that 

$$
f\in \prod\limits_{i<\kappa}\lambda_i \Rightarrow 
(\exists \gamma_{f}<\lambda) (\{i<\kappa:f(i)< \gamma_f\}\in D).
$$

\par \noindent \Then\ $D$ is $(\kappa,g)$-weakly normal (hence $(\kappa,\theta)$-weakly normal). 
\end{claim}

\par \noindent \emph{Proof}. \newline 
We will show that $D$ is $(\kappa,g)$-weakly 
normal. Let $h:\kappa\rightarrow \theta$ be any 
$(\kappa,g)$-regressive function.

\par \noindent For every $i<\kappa$ define $f(i)=\mu_{h(i)}+1$. 
Clearly $f\in \prod\limits_{i<\kappa}\lambda_i$
since $\lambda_i=\mu_{g(i)}$ and $\mu_{h(i)}<
\mu_{g(i)}$ for every $i<\kappa$. 
Let $\gamma_f<\lambda$ be \st\ 
$\{i<\kappa:f(i)<\gamma_f\}\in D$. Define $j_h$ 
to be the first ordinal \st\ 
$\mu_{j_h}>\gamma_f$. By that, we have
$\{i<\kappa:h(i)<j_h\}\in D$, so we are done. \newline 

\hfill \qedref{0.8}

\medskip

We have defined some notions of normality, 
for ultrafilters. The other side of the coin is 
regular ultrafilters. A good source to this subject is \cite{MR0409165}. Let us start with the definition:

\begin{definition}
\label{0.9}
Regular ultrafilters. \newline 
Let $D$ be an ultrafilter on $\kappa, \alpha\leq \kappa$.
\begin{enumerate}
\item[(a)] $D$ is $\alpha$-regular if there exists 
$E\subseteq D, |E|=\alpha$, and for every 
$i<\kappa$ we have $|\{e\in E:i\in e\}|<\aleph_0$
\item[(b)] $D$ is regular, when $\alpha=\kappa$
\end{enumerate}
\end{definition}

Notice that every ultrafilter is $\alpha$-regular 
for any $\alpha<\aleph_0$, so the definition is 
interesting only when $\alpha$ is an infinite cardinal. 
But even in the first infinite cardinal, i.e. 
$\alpha=\aleph_0$, we have a useful result for 
our needs. 

\begin{claim}
\label{0.10}
An ultrafilter $D$ on $\kappa$ is $\aleph_0$-regular iff 
it is not $\aleph_1$-complete.
\end{claim}

\par \noindent \emph{Proof}. \newline 
If $D$ is $\aleph_0$-regular, let $E\subseteq D$ 
be an evidence.
Every $i<\kappa$ belongs to a finite subset of $E$, 
and $|E|=\aleph_0$, so $i\notin \bigcap E$ for any 
$i<\kappa$. In other words, 
$\bigcap E=\emptyset \notin D$, so $D$ is not 
$\aleph_1$-complete.

If $D$ is not $\aleph_1$-complete, we can find a 
countable $E\subseteq D$, \st\ 
$\bigcap E\notin D$. Leaning on the fact that $D$ 
is an ultrafilter, we can define a countable $E'$ 
which stands in the demands of the $\aleph_0$ 
regularity. \newline 

\hfill \qedref{0.10}

\par \noindent We state the following well-known results, without a proof:

\begin{theorem}
\label{0.11}
Let $\kappa$ be the first cardinal \st\ we have a non-principal 
$\aleph_1$-complete ultrafilter on it. 
\Then\ $\kappa$ is a measurable cardinal. \newline 
\end{theorem}

\hfill \qedref{0.11}

\begin{theorem}
\label{solovay}
Suppose $\mu$ is a compact cardinal, $\chi = \cf(\chi) \geq \mu$, and $\theta < \mu$. \Then\ $\chi^\theta = \chi$.
\end{theorem}

\hfill \qedref{solovay}

\par \noindent The proof of these theorems can be found in \cite{MR1321144}.

We conclude this section with an important cardinal 
arithmetic result, for $\aleph_0$-regular 
ultrafilters (the proof can be found in \cite{MR0409165}):

\begin{claim}
\label{0.12}
Let $A$ be an infinite set, $D$ an $\aleph_0$-regular
ultrafilter on $\tau$. \newline 
\Then\ $|\prod\limits_{\tau} A /D|\geq |A|^{\aleph_0}$. \newline 
\end{claim}

\hfill \qedref{0.12}

We thank the referee for the excellent work, which was much deeper than just simple proofreading.

\newpage

\section{weak normality and low cardinality}

The title of this section is not just a rhyme. 
It captures mathematical information. For showing 
this, let us start with the simple direction.

\begin{proposition}
\label{1.1}
Assume $D$ is a $(\kappa,\theta)$-weakly normal 
ultrafilter on $\kappa$,
$\cf (\lambda)=\theta$ and $\lambda$ is $\kappa$-strong. \newline 
\Then\ we can find a \seq\ of cardinals
$\bar \lambda=\langle \lambda_i:i<\kappa\rangle$
\st\ $\lambda={\rm lim}_D (\bar \lambda)$ and
$|\prod\limits_{i<\kappa} \lambda_i/D|=\lambda$.
\end{proposition}

\par \noindent \emph{Proof}. \newline 
First, we choose our sequence. Let $\bar \mu=\langle \mu_j:j<\theta\rangle$ be a 
continuous \incr\ \seq\ of cardinals, with limit $\lambda$.
Let $g:\kappa\rightarrow \theta$ be a witness to the 
$(\kappa,\theta)$-weak normality of $D$. Define 
$\lambda_i=\mu_{g(i)}$, for any $i<\kappa$.
By \ref{0.5} we know that $\lambda={\rm lim}_D 
(\langle\lambda_i:i<\kappa\rangle)$.

\par \noindent Now, we must prove two inequalities:
\begin{enumerate}
\item[(a)] $\lambda\leq |\prod\limits_{i<\kappa}
\lambda_i/D|$

\par \noindent By \ref{0.6} and the fact that $\lambda={\rm lim}_D 
(\langle\lambda_i:i<\kappa\rangle)$, we conclude that $\lambda\leq |\prod\limits_{i<\kappa} \lambda_i/D|$.

\item[(b)] $|\prod\limits_{i<\kappa} \lambda_i
/D|\leq \lambda$

\par \noindent Observe that for every $f\in \prod\limits_{i<\kappa}
\lambda_i$ we can find $\gamma_f<\lambda$, \st\ 
$\{i<\kappa:f(i)\leq \gamma_f\}\in D$.

\par \noindent Why? Well, $f\in \prod\limits_{i<\kappa} \lambda_i
=\prod\limits_{i<\kappa} \mu_{g(i)}$. 
Define $f^*:\kappa\rightarrow \theta$ in 
the following way: for every $i<\kappa$ let $f^*(i)$ 
be the first ordinal $j$ \st\ $f(i)<\mu_j$. 
$f^*$ is $(\kappa,g)$-regressive (truely, we have 
$f^*(i)\leq g(i)$, but the difference between $\leq$ 
and $<$ is unimportant here). By the $(\kappa,g)$-weak normality assumption, one can find $j<\theta$ \st\ the set $\{i<\kappa:f^*(i)<j\}\in D$. That means 
also that the set $\{i<\kappa:f(i)<\mu_j\}$ 
belongs to $D$, so choose $\gamma_f =\mu_j$ and 
the assertion follows.

\par \noindent For $\gamma < \lambda$ let $\mathcal{F}_\gamma$ be the set $\{ f\in 
\prod\limits_{i<\kappa}\lambda_i$, 
and $\{i<\kappa:f(i)<\gamma\}\in D\}$. Now, we have: \newline 
$
|\prod\limits_{i<\kappa}\lambda_i/D|=
|\{f/D:f\in \prod\limits_{i<\kappa} \lambda_i\}|\leq
|\bigcup\limits_{\gamma<\lambda}f/D:f\in 
\mathcal{F}_\gamma|\leq 
\sum\limits_{\gamma<\lambda} |\{f/D:f\in 
\mathcal{F}_\gamma| \leq \sum\limits_{\gamma<\lambda} |\gamma|^\kappa
\leq \lambda\times \lambda=\lambda$.
\end{enumerate}

\hfill \qedref{1.1}

\medskip

One remark about proposition \ref{1.1}. 
We took an infinite $\lambda$ \st\ 
$\lambda$ is $\kappa$-strong.
Clearly, that assumption is vital, since 
$\lambda^\kappa>\lambda$ in our case. 
So under that necessary restriction on $\lambda$,
all we need for the low cardinality of the product is 
the $(\kappa,\theta)$-weak normality of $D$.

\par \noindent We turn now to the opposite direction:

\begin{theorem}
\label{1.2}
Assume 
\begin{enumerate}
\item[(a)] $\theta=\cf(\lambda)\leq \kappa<\lambda$
\item[(b)] $\langle\lambda_i:i<\kappa\rangle$ is a sequence of cardinals
\item[(c)] ${\rm lim}_D  (\langle\lambda_i:i<\kappa\rangle)=\lambda$
\item[(d)] $\lambda^\kappa_i<\lambda$ for every $i<\kappa$
\item[(e)] $D$ is an ultrafilter on $\kappa$
\item[(f)] $D$ is not closed to descending sequences of length $\theta$ (e.g., $D$ is not $\aleph_1$-complete)
\item[(g)] $|\prod\limits_{i<\kappa} \lambda_i /D|=\lambda$
\end{enumerate}
\Then\ $D$ is $(\kappa,\theta)$-weakly normal.
\end{theorem}

\par \noindent \emph{Proof}. \newline 
Let $\bar f=\langle f_\alpha:\alpha<\lambda\rangle$ 
be a set of representatives to $\prod\limits_{i<\kappa}
\lambda_i/D$. Denote $\kappa_i=\cf(\lambda_i)$,
for every $i<\kappa$.
\begin{enumerate}
\item[$(*)_0$] ${\rm lim}_D (\langle\kappa_i:i<\kappa\rangle)
<\lambda$.

\par \noindent [Why? If ${\rm lim}_D (\langle\kappa_i:i<\kappa\rangle)\geq \lambda$ then $\beta < \lambda \Rightarrow \{ i < \kappa : \beta < \kappa_i \} \in D$, so $(\prod\limits_{i<\kappa} \kappa_i,
\leq_D)$ is $\lambda$-directed, by \ref{0.6} and \ref{0.7}(b).
Since $\lambda$ is singular, $(\prod\limits_{i<\kappa}
\kappa_i,\leq_D)$ is even $\lambda^+$-directed, so by 
\ref{0.7}(a) $|\prod\limits_{i<\kappa}
\kappa_i/D|\geq \lambda^+$ and consequently $|\prod\limits_{i<\kappa}
\lambda_i/D|\geq \lambda^+$ since $\kappa_i\leq\lambda_i$ for every $i<\kappa$, contradicting assumption (f) here].

\par \noindent It follows from $(*)_0$ that $\kappa_i=\cf(\lambda_i)<\lambda_i$
for a set of $i$'s which belongs to $D$. \Wlog, we can
assume that:
\item[$(*)_1$] $\cf(\lambda_i)<\lambda_i$, for every 
$i<\kappa$.

\par \noindent For every $i<\kappa$, choose 
$\langle \lambda_{i,\epsilon} : \epsilon < \kappa_i \rangle$ such that:
\begin{enumerate}
\item[(i)] $\kappa<\lambda_{i,\epsilon}=
\cf(\lambda_{i,\epsilon})<\lambda_i$
\item[(ii)] $\sum\limits_{\epsilon<\kappa_i}
\lambda_{i,\epsilon}=\lambda_i$
\item[(iii)] $\kappa_{i_1}=\kappa_{i_2}
\Rightarrow \lambda_{i_1,\epsilon}=\lambda_{i_2,\epsilon}$
for every $\epsilon<\kappa_{i_1}=\kappa_{i_2}$
\end{enumerate}
\item[$(*)_2$] There is no $h\in \prod\limits_{i<\kappa}
\kappa_i$, \st\  ${\rm lim}_D (\langle\lambda_{i,h(i)}:i<\kappa\rangle)
=\lambda$. 

\par \noindent [Why? exactly like $(*)_0$, upon replacing $\kappa_i$ by $\lambda^+_{i, h(i)}$] 

\par \noindent Let $\bar \mu=\langle \mu_j:j<\theta\rangle$ be an \incr\ 
\cont\ \seq\ of singular cardinals, with limit $\lambda$. 
Notice that $\theta>\aleph_0$ here, by (f) and (g), hence such $\bar{\mu}$ exists. \newline 
We claim that for $D$-many $i$'s we have $\lambda_i\in 
\{\mu_j:j<\theta\}$. Otherwise, define 
$\zeta(i)={\rm sup}\{j:\mu_j<\lambda_i\}$ for $i<\kappa$. 
Since $\bar{\mu}$ is continuous, we will get $\mu_{\zeta(i)}<\lambda_i$ for $D$-many $i$'s, so easily one can create $h\in 
\prod\limits_{i<\kappa}
\kappa_i$ \st\ $\mu_{\zeta(i)}<\lambda_{i,h(i)}<\lambda_i$.
Clearly, we have ${\rm lim}_D(\langle \lambda_{i,h(i)}
:i<\kappa\rangle)=\lambda$, contradicting $(*)_2$.

\par \noindent So, without loss of generality:
\item[$(*)_3$] $\lambda_i\in \{\mu_j:j<\theta\}$, 
for every $i<\kappa$. 

\par \noindent For each $i<\kappa$, let $g(i)$ be the first ordinal 
$j<\theta$ \st\ $\lambda_i=\mu_j$.
We will show (in $(*)_4$ below) that $g$ is a witness to the 
$(\kappa,\theta)$-weak normality of $D$.

\item[$(*)_4$]
For every $f\in \prod\limits_{i<\kappa}\lambda_i$,
there is $\gamma_f<\lambda$, such that: 
$$
\{i<\kappa:f(i)<\gamma_f\}\in D
$$
\end{enumerate}

\par \noindent For every $i<\kappa$, define 
$P_i=\{\lambda_{i,\epsilon}:\epsilon<\kappa_i\}$.
By the choice of the $\lambda_{i,\epsilon}$-s,
$\prod\limits_{i<\kappa}P_i/D$ is unbounded in 
$\prod\limits_{i<\kappa}\lambda_i/D$. 
Observe that tcf$(\prod\limits_{i<\kappa} P_i
/D)=\cf(\lambda)=\theta$, since 
$|\prod\limits_{i<\kappa} \lambda_i/D|=\lambda$. 
Consequently, tcf$(\prod\limits_{i<\kappa} \kappa_i
/D)=\theta$, since ${\rm otp} (P_i)=
\kappa_i$ for every $i<\kappa$.

Now, let $f\in \prod\limits_{i<\kappa} \lambda_i$
be any function. $\prod\limits_{i<\kappa}
P_i/D$ is unbounded in $\prod\limits_{i<\kappa}
\lambda_i/D$, so we can find 
$m\in \prod\limits_{i<\kappa} P_i/D$ \st\ $f<_D m$.
By $(*)_2$ and the observation above, 
$\gamma:={\rm lim}_D (\langle m(i):i<\kappa\rangle)
<\lambda$. Choose $\gamma_f=\gamma$, and the proof of $(*)_4$ is complete.

Now we can finish the proof of the theorem. 
Just notice that claim \ref{0.8} asserts,
under $(*)_4$, that $D$ is 
$(\kappa,\theta)$-weakly normal 
(with respect to the function $g$, which is defined above).

\hfill \qedref{1.2}

We conclude this section with the case of singular cardinals with countable cofinality. One of the early results about the continuum hypothesis, 
much before the Cohen era and even before G\"odel, asserts that $2^{\aleph_0}\neq \aleph_\omega$. More generally, if $\cf(\lambda)=\aleph_0$, then $\lambda$ can not realize the continuum (the result belongs to K\"onig, and appears in \cite{MR1511338}).

One of the metamathematical ideas of the pcf theory, suggests to replace the questions of $2^\theta$ by questions of products of cardinals, modulo an ultrafilter. We would like to phrase a similar result
about singular $\lambda$-s with countable cofinality, 
this time in the light of the pcf. This result is the content of corollary \ref{morereferee} below.

\begin{proposition}
\label{rreferee}
$(\kappa, \aleph_0)$-weak normality. \newline 
For any cardinal $\kappa$ there is no $(\kappa, \aleph_0)$-weakly normal ultrafilter on $\kappa$.
\end{proposition}

\par \noindent \emph{Proof}. \newline 
Suppose that $D$ is an ultrafilter on $\kappa$, and $g : \kappa \rightarrow \aleph_0$ satisfies condition $(i)$ of definition \ref{0.3} $(a)$. It means that $\{i\in\kappa:g(i)>j\}\in D$ for every $j\in\omega$.
Let $f : \kappa \rightarrow \aleph_0$ be defined by $f(i)=g(i)-1$ (and if $g(i)=0$ then $f(i)=g(i)$) for all $i<\kappa$. 

Then for every $j \in \omega \setminus \{0\}$ we have $\{i \in \kappa : f(i) < j\} = \{i \in \kappa : g(i) \leq j\} = \aleph_0 \setminus \{i \in \kappa : g(i) > j\} \notin D$, by $(i)$.

\hfill \qedref{rreferee}

\begin{corollary}
\label{morereferee}
Assume
\begin{enumerate}
\item [$(a)$] $\aleph_0 = \cf(\lambda) \leq \kappa < \lambda$
\item [$(b)$] $\bar{\lambda} = \langle \lambda_i : i < \kappa \rangle$ is a sequence of cardinals
\item [$(c)$] ${\rm lim}_D (\bar{\lambda}) = \lambda$
\item [$(d)$] $\lambda$ is $\kappa$-strong
\item [$(e)$] $D$ is an ultrafilter on $\kappa$
\item[(f)] $D$ is not closed to descending sequences of length $\theta$ (e.g., $D$ is not $\aleph_1$-complete)
\end{enumerate}
\Then\ $|\prod \limits_{i<\kappa} \lambda_i / D| \neq \lambda$
\end{corollary}

\hfill \qedref{morereferee}

\begin{remark}
\label{aleph1}
Measurability and weak normality.
\begin{enumerate}
\item [$(\aleph)$] If $D$ is closed under descending sequences of length $\cf(\lambda)$, then theorem \ref{1.2} and the former corollary need not to be true. In particular, if $\kappa$ is a measurable cardinal and $\lambda>\kappa$, $\cf(\lambda)=\aleph_0$, then $\lambda$ can be realized as $|\prod \limits_{i<\kappa} \lambda_i / D|$. Nevertheless, $D$ is not $(\kappa,\aleph_0)$-weakly normal (see \ref{rreferee}).
\item [$(\beth)$] The situation is different for $(\kappa, \theta)$-weakly normal ultrafilters when $\theta > \aleph_0$. For example, it is consistent to have a weakly normal (uniform) ultrafilter on $\aleph_1$, see \cite{MR942519} and the history there, and see also \cite{MR1713438}.
\end{enumerate}
\end{remark}

\newpage

\section{Applications to Boolean algebras}

\par \noindent We turn now to the field of Boolean Algebras:

\begin{definition}
\label{3.1}
$\Depth$ and $\Depth^+$. \newline 
Let $\bool$ be a Boolean Algebra. 
\begin{enumerate}
\item[(a)] $\Depth(\bool):={\rm sup}
\{\theta:$ there exists $A\subseteq \bool,|A|=\theta$,\newline
\qquad \qquad $A$ is well ordered by $<_\bool\}$
\item[(b)] ${\rm Depth}^+ (\bool):={\rm sup}\{\theta^+$:
there exists $A\subseteq \bool, |A|=\theta$,\newline
\qquad\qquad $A$ is well-ordered by $<_\bool\}$
\end{enumerate}
\end{definition}

\par \noindent Monk raised the following question:

\begin{question}
\label{3.2}
Let $\langle \bool_i:i<\theta\rangle$ be a \seq\ of Boolean Algebras,
$D$ an ultrafilter on $\theta, \bool=\prod\limits_{i<\theta}
\bool_i/D$. 
Can we have, in ZFC, an example of ${\rm Depth}(\bool)
> \prod\limits_{i<\theta} {\rm Depth} (\bool_i)/D$?
\end{question}

We try to find a necessary condition for such an example
above a compact cardinal. We start with the following claim, 
from \cite{MR2217966}:

\begin{claim}
\label{3.3A}
Assume 
\begin{enumerate}
\item[(a)] $\theta<\mu\leq \lambda$
\item[(b)] $\mu$ is a compact cardinal
\item[(c)] $\lambda=\cf(\lambda),D$ is an ultrafilter 
on $\theta$
\item[(d)] $(\forall \alpha<\lambda) 
(|\alpha|^\theta <\lambda)$
\item[(e)] ${\rm Depth}^+ (\bool_i)\leq \lambda$ 
for every $i<\theta$
\end{enumerate}
\Then\ ${\rm Depth}^+ (\bool)\leq \lambda$.
\end{claim}

\hfill \qedref{3.3A}

As a simple conclusion, we can derive our necessary 
condition in terms of cardinal arithmetic:

\begin{conclusion}
\label{3.4}
Assume 
\begin{enumerate}
\item[(a)] $\theta<\mu<\lambda$
\item[(b)] $\mu$ is a compact cardinal 
\item[(c)] ${\rm Depth} (\bool_i)\leq \lambda$, 
for every $i<\theta$
\item[(d)] $D$ is a uniform ultrafilter on $\theta$
\item[(e)] ${\rm Depth}(\bool) > \lambda + \prod\limits_{i<\theta} {\rm Depth} (\bool_i)/D$
\end{enumerate}
\Then\ ${\rm Depth}(\bool)=\lambda^+$ and 
$|\prod\limits_{i<\theta} {\rm Depth} (\bool_i)/D|=\lambda$.
\end{conclusion}

\par \noindent \emph{Proof}. \newline 
By (c) we know that $\Depth^+ (\bool_i)\leq \lambda^+$, 
so clearly $\Depth^+(\bool_i)\leq \lambda^{++}$ for every 
$i<\theta$. Now, $\lambda^{++}$ stands in the demands of 
claim \ref{3.3A} (remember that $\mu$ is compact, 
so $(\lambda^+)^\theta=\lambda^+$ by Solovay's theorem). 
Hence $\Depth^+(\bool) \leq \lambda^{++}$, and consequently 
$\Depth(\bool)\leq \lambda^+$. 

By (e), $|\prod\limits_{i<\theta} \Depth(\bool_i)/D|$ 
is strictly less than $\Depth(\bool)$. 
But we deal here with a case of $\prod\limits_{i<\theta}
\Depth(\bool_i)/D\geq \lambda$, so the only 
possibility is $\Depth(\bool)=\lambda^+$ and 
$|\prod\limits_{i<\theta}\Depth(\bool_i)/D|=\lambda$.

\hfill \qedref{3.4}

\medskip

We focus, from now on, in the case of a singular $\lambda$
with cofinality $\aleph_0$. In general, it seems that 
those cardinals behave in a unique way around questions 
of Depth. The following Theorem shows that there is a limitation on examples like \ref{3.4}, for a singular $\lambda$ with countable cofinality:

\begin{theorem}
\label{3.5}
Assume 
\begin{enumerate}
\item[(a)] $\kappa<\mu<\lambda,\mu$ is a compact cardinal
\item[(b)] $\kappa$ is the first measurable cardinal, $\theta<\kappa$
\item[(c)] $\lambda$ is a singular cardinal, 
$\cf(\lambda)=\aleph_0$
\item[(d)] $\langle \bool_i:i<\theta\rangle$ is a \seq\ 
of Boolean Algebras
\item[(e)] $D$ is a uniform ultrafilter on $\theta$
\item[(f)] $\Depth(\bool_i) \leq \lambda$, for every $i<\theta$
\end{enumerate}
\Then\ $\Depth(\bool)\leq \lambda$.
\end{theorem}

\par \noindent \emph{Proof}. \newline 
Assume toward contradiction, that 
$\Depth(\bool)>\lambda$. Due to \ref{3.4}, we have 
an example of $|\prod\limits_{i<\theta} \Depth(\bool_i)/D| < \Depth(\bool)$ above a compact $\mu$, so by 
virtue of conclusion \ref{3.4} we must have 
$|\prod\limits_{i<\theta} \Depth(\bool_i)/D|=\lambda$. Theorem \ref{1.2} implies, under this consideration, that $D$ is $(\theta, \aleph_0)$-weakly normal. But this is impossible, as shown in \ref{rreferee}.

\hfill \qedref{3.5}

\begin{remark}
\label{magidor}
Consistency results.
\begin{enumerate}
\item [(a)] By \cite{MR0429566} it is consistent that the first compact is the first measurable. Consequently, there is no example of $|\prod\limits_{i<\theta} \Depth(\bool_i)/D| < \Depth(\bool)$ for singular $\lambda$-s with countable cofinality above the first measurable cardinal, in ZFC.
\item [(b)] By \cite{temp911}, if $\bf {V} = \bf {L}$ there is no example as above. This paper gives (part of) the picture under large cardinals assumptions.
\end{enumerate}
\end{remark}

\newpage 

\bibliographystyle{amsplain}
\bibliography{arlist}

\end{document}